\documentclass[12pt]{amsart}
\usepackage{amssymb,amsbsy,amsthm,amscd}
\usepackage[mathcal]{eucal}

\usepackage[all]{xy}

\theoremstyle{plain}
\newtheorem{thm}{Theorem}[section]
\newtheorem{lem}[thm]{Lemma}
\newtheorem{cor}[thm]{Corollary}
\newtheorem{prop}[thm]{Proposition}

\numberwithin{equation}{section}
\newcommand{\N}{\ensuremath{\mathbb N}}
\newcommand{\Z}{\ensuremath{\mathbb Z}}
\newcommand{\F}{\ensuremath{\mathbb F}}

\DeclareMathOperator{\GL}{GL}
\DeclareMathOperator{\SL}{SL}
\DeclareMathOperator{\ord}{ord}
\DeclareMathOperator{\Eig}{Eig}
\DeclareMathOperator{\Pasc}{Pasc}
\DeclareMathOperator{\tr}{tr}
\DeclareMathOperator{\RC}{RC}
\DeclareMathOperator{\diag}{diag}

\begin{document}
\date\today
\title{Regular Circulant Matrices}
\thanks{The author would like to thank F.~Grunewald and B.~Klopsch for helpful discussions. The author was supported by a Thomas Holloway Scholarship of Royal Holloway College, University of London.}
\author{Daniel Appel}
\address{Department of Mathematics, Royal Holloway College, University of London, Egham, Surrey, TW20 0EX,
United Kingdom.}
\curraddr{Mathematisches Institut der Heinrich- Heine- Universit\"at\\ 40225 D\"usseldorf\\ Germany.}
\email{Daniel.Appel@uni-duesseldorf.de}

\maketitle
\pagestyle{plain}
\begin{abstract}
We consider the groups $\RC_n(\F_{p^t})$ and $\RC_n(\Z / a \Z)$ of regular circulant $(n \times n)$-matrices over $\F_{p^t}$ and $\Z / a\Z$, respectively, where $p$ is a prime and $t,n,a \in \N$. In both cases we present a formula for the order of that group. We also make a first step towards finding the algebraic structure of these groups.
\end{abstract}

\section{Introduction and Main Results}
We consider regular circulant matrices over finite fields and integer residue class rings. In general, a matrix $\begin{pmatrix} a_{i,j} \end{pmatrix}_{1 \leq i,j \leq n}$ is called \emph{circulant}, if, for all $1 \leq i,j \leq n$, we have $a_{i,j+1} = a_{i-1,j}$ where the indices have to be read modulo $n$. Hence, a circulant matrix is completely determined by any of its columns (respectively rows) and each column (respectively row) can be obtained from the previous one by a cyclic permutation.

It is commonly known that the product of two circulant matrices is again circulant and so is the inverse of a regular circulant matrix. Therefore we may consider the groups $\RC_n(\F_{p^t})$ and $\RC_n(\Z / a \Z)$ of regular circulant $(n \times n)$-matrices over $\F_{p^t}$ and $\Z / a\Z$, respectively, where $p$ is a prime and $t,n,a \in \N$. 

Let us set up some notation that we need to state our results. Suppose that $\F_{p^s} \geq \F_{p^t}$ is a field extension. Then we write $F^t : \F_{p^s} \rightarrow \F_{p^s}$, $x \mapsto x^{p^t}$ for the relative Frobenius Homomorphism. To vectors and matrices we apply $F^t$ componentwise. For a vector $w = \begin{pmatrix} w_i \end{pmatrix} \in \F_{p^s}^{k}$, we define an upper triangular matrix $T(w) :=  \begin{pmatrix} w_{k-j+i}\end{pmatrix}_{1\leq i,j \leq k}$ where $w_\lambda := 0$ for $\lambda \leq 0$. Moreover, if $d,k \in \N$ with $\gcd(d,k)=1$, we write $\ord_d(k)$ for the order of $k$ in $(\Z/d\Z)^*$. Finally, by $\phi$ we denote the Euler function.

\begin{thm}\label{main1}
Let $p$ be a prime and $n =  mp^r \in \N$ with $p \nmid m$.  
\begin{itemize}
\item[(i)] The order of $\RC_{n}(\F_{p^t})$ is
$$\prod_{d \mid m} (p^{t\cdot \ord_d(p^t)} - 1)^{\phi(d)/\ord_d(p^t)} \cdot p^{t \cdot (p^r - 1)\cdot \phi(d)}.$$
\item[(ii)] Let $\F_{p^s} \geq \F_{p^t}$ be a field extension such that $\F_{p^s}$ contains the $m$-th roots of unity. Write the permutation $\sigma : \Z/m\Z \rightarrow \Z/m\Z$, $x \mapsto p^{-t}x$ as $\sigma = \sigma_1 \cdots \sigma_l$ with disjoint cycles $\sigma_k = (s_{k,1},\dots,s_{k,m_k})$ of length $m_k$. Then, as a subgroup of $\GL_n(\F_{p^s})$, the group $\RC_n(\F_{p^t})$ is conjugate to the group consisting of the matrices
$$
\begin{pmatrix} T(v_{(1)}) \\ & \ddots \\ & & T(v_{(m)}) \end{pmatrix}
$$
with $v_{(i)} = \begin{pmatrix} v_{i,j} \end{pmatrix} \in \F_{p^s}^{p^r}$  satisfying
\begin{align*}
& v_{1,p^r}, v_{2,p^r}, \dots v_{m,p^r} \not= 0 \mbox{ and}\\
& v_{(s_{k,1})} \in \F_{p^{tm_k}}^{p^r},\  v_{(s_{k,i})} = F^{(i-1)t}(v_{(s_{k,1})})
\end{align*}
for $2 \leq i \leq m_k$, $1 \leq k \leq l$.
\end{itemize}
\end{thm}

\begin{thm}\label{main2}
Let $a, n\in \N$ and $a = \prod_{p \mid a} p^{t_p}$ be the prime factorization of $a$. Moreover, for every prime divisor $p$ of $a$, define $r_p, m_p \in \N$ by $n = p^{r_p}\cdot m_p$ with $p \nmid m_p$. Then the order of $\RC_n(\Z / a \Z)$ is
$$\prod_{p \mid a} \left( p^{n(t_p - 1)} \prod_{d \mid m_p} (p^{\ord_d(p)} - 1)^{\phi(d)/\ord_d(p)} \cdot p^{(p^{r_p} - 1)\cdot \phi(d)} \right).$$
\end{thm}

We use very explicit computations to prove our results. An alternative approach would be to consider the algebra isomorphism between $\RC_n(\F_{p^t})$ and $\F_{p^t}[x]/(x^n-1)$. This approach has been used to determine the number orthogonal $(n \times n)$-matrices over $\F_{p^t}$. See for example \cite{JBG} and the references therein for further details.

From Theorem \ref{main1} we also obtain
\begin{cor}
Let $p$ be a prime and $n,t \in \N$. The group $\RC_n(\F_{p^t}) \cap \SL_n(\F_{p^t})$ has index $p^t-1$ in $\RC_n(\F_{p^t})$. In particular, there are exactly $|\RC_n(\F_p^t)| \cdot (p^t-1)^{-1} $ circulant $(n \times n)$-matrices of determinant $1$ over~$\F_{p^t}$
\end{cor}

\begin{proof}
We show that $\det : \RC_n(\F_{p^t}) \rightarrow \F_{p^t}^*$ is onto. Since $\RC_n(\F_{p^t}) \cap \SL_n(\F_{p^t})$ is the kernel of this homomorphism, this proves the claim.

Let $X$ be a matrix that conjugates $\RC_n(\F_{p^t})$ to the group given in part (ii) of Theorem \ref{main1}. It suffices to show that $\det:X \RC_n(\F_{p^t}) X^{-1} \rightarrow \F_{p^t}^*$ is onto. Observe that $\sigma$ has the fixed point $m \in \Z / m\Z$. Therefore the vector $v_{(m)} \in \F_{p^t}^{p^r}$ can be chosen independently from $v_{(1)}, \dots, v_{(m-1)}$, subject only to the restriction that $v_{m,p^r} \not= 0$. Note that $\det(T(v_{(m)})) = v_{(m,p^r)}^{p^r}$. Since $\gcd(p^r,p^t-1)=1$, the map $\F_{p^t} \rightarrow\F_{p^t}$, $x \mapsto x^{p^r}$ is a bijection. Hence, every element of $\F_{p^t}^*$ can be obtained as $\det(T(v_{(m)}))$ with suitable $v_{(m)}$. This implies that every element of $\F_{p^t}^*$ can be obtained as the determinant of a matrix as in part (ii) of Theorem \ref{main1}.
\end{proof}

\section{Preliminaries}\label{prel}
\subsection{Permutation induced by the Frobenius Homomorphism}

Let $p$ be a prime and $m,t \in \N$ such that $p \nmid m$. Moreover, let $\F_{p^s} \geq \F_{p^t}$ be a field extension of $\F_{p^t}$ that contains the $m$-th roots of unity and let 
$$F^t: \F_{p^s} \longrightarrow \F_{p^s},\quad x \longmapsto x^{p^t}$$
be the relative Frobenius Homomorphism. Observe that $F^t$ induces a permutation of the $m$-th roots of unity, which is described by
\begin{equation}\label{sigma} \sigma : \Z / m\Z \longrightarrow \Z / m\Z, \quad a \longmapsto p^t a.\end{equation}
We consider the cycle structure of this permutation. To this end, we first set up some notation.

For $d \in \N$ with $p \nmid d$ we write $\ord_d(p^t)$ for the order of $p^t$ in $(\Z / d\Z)^*$. By $\langle p^t \rangle \leq (\Z/m\Z)^*$ we denote the subgroup of $(\Z/m\Z)^*$ generated by $p^t$. Moreover, for $d \mid m$ we set $V_d := d \cdot \langle p^t \rangle \subseteq d \cdot (\Z / m\Z)^*$. In order to describe the permutation $\sigma$, we consider the action of $\langle p^t \rangle$ on $\Z / m\Z$ by multiplication.

Observe that
$$\Z / m\Z = \biguplus_{d \mid m} d \cdot (\Z/m\Z)^*.$$
Moreover, each set $d \cdot (\Z/m\Z)^*$ is a disjoint union of $|d \cdot (\Z/m\Z)^*| / |V_d|$ sets of the form $\varepsilon \cdot V_d$ with certain $\varepsilon \in (\Z/m\Z)^*$. In particular, we have
$$\Z /m\Z = \biguplus_{d \mid m} \biguplus_{\varepsilon} \varepsilon \cdot V_d.$$
Clearly the sets $\varepsilon \cdot V_d = \{ \varepsilon d p^{tk} \mid k \in \Z \}$ are invariant under the action of $\langle p^t \rangle$ and the action of $\langle p^t \rangle$ on each set $\varepsilon \cdot V_d$ is given by a cycle of length $|V_d|$. Hence, the action of $\langle p^t \rangle$ on $d \cdot (\Z /m\Z)^*$ is described by a product of $|d\cdot(\Z/m\Z)^*|/|V_d|$ disjoint cycles of length $|V_d|$, that is, it has cycle structure $(|V_d|)^{|d \cdot(\Z/m\Z)^*|/|V_d|}$. This shows that the permutation $\sigma$ in (\ref{sigma}) has cycle structure
$$\prod_{d \mid m} (|V_d|)^{|d \cdot (\Z/m\Z)^*|/|V_d|}.$$

Now observe that for $d \mid m$ we have $|d \cdot (\Z / m\Z)^*| = \phi(m/d)$, where $\phi$ denotes the Euler function. One also easily verifies that $|V_d| = \ord_{m/d}(p^t)$. Noting that, as $d$ runs through all divisors of $m$, so does $m/d$, we thus obtain
\begin{prop}\label{cycle}
The permutation $\sigma$ in (\ref{sigma}) has cycle structure
$$ \prod_{d\mid m} (\ord_d(p^t))^{\phi(d) / \ord_d(p^t)}.$$
\end{prop}

\subsection{The Kronecker Product}
Let us briefly recall the notion and some basic properties of the Kronecker product of matrices. For details we refer to \cite{WS}.

Given an $(m \times n)$-matrix $A = \begin{pmatrix} a_{ij} \end{pmatrix}$ and an $(r \times s)$-matrix $B$, the Kronecker product $A \otimes B$ of $A$ and $B$ is the $(mr \times ns)$-matrix $A \otimes B := \begin{pmatrix}  a_{ij} B \end{pmatrix}$. By $A^{\otimes k}$ we denote the $k$-fold Kronecker product of $A$ with itself. Two standard result are
$$ (A \otimes B)^{\tr} = A^{\tr} \otimes B^{\tr} \quad \mbox{and} \quad (A \otimes B)^{-1} = A^{-1} \otimes B^{-1}.$$
Observe that these imply that 
$$(A^{\otimes k})^{\tr} = (A^{\tr})^{\otimes k} \quad \mbox{and} \quad (A^{\otimes k})^{-1} = (A^{-1})^{\otimes k}.$$
Now let $A$, $B$, $C$, $D$ all be matrices of the same dimension. Then another standard result says $AC \otimes BD = (A \otimes B) (C \otimes D)$. If $A$ is regular, one easily obtains from this that
$$
(ABA^{-1})^{\otimes k} = A^{\otimes k} B^{\otimes k} (A^{-1})^{\otimes k}.
$$

\subsection{Pascal Matrices modulo primes}
For $n \in \N$ we let $\Pasc_n(\Z)$ be the Pascal matrix of size $n \times n$ over $\Z$, that is,
$$\Pasc_n(\Z) = \begin{pmatrix} {i-1 \choose j-1} \end{pmatrix}_{i,j}. $$
Let $p$ be a prime. Then, through the natural projection $\Z \rightarrow \F_p$, we can also define the Pascal matrix $\Pasc_n(\F_p)$ over $\F_p$.

It is well known that
\begin{equation}\label{pasc}
\Pasc_n(\Z)^{-1} = \begin{pmatrix} (-1)^{i+j} {i-1 \choose j-1} \end{pmatrix}_{i,j}.
\end{equation}
See for example \cite{CV} for this result. As we shall see, the inverse of $\Pasc_{p^r}(\F_p)$, with $r \in \N$, can also be described by a different nice formula.
\begin{lem}
Let $p$ be a prime. Then ${p-j \choose p-i} \equiv (-1)^{i+j} {i-1 \choose j-1} \mod p$ for $1 \leq i,j \leq p$.
\end{lem}

\begin{proof}
For $j > i$ both terms are zero. So we may assume that $i \geq j$. Set $k:=i-j$ so that $i=j+k$ and $0 \leq k \leq p-1$. Observing that $\gcd(p,k!) = 1$, we thus find that our claim is equivalent to
\begin{align*}
&		&			 					&\quad {p-j \choose p-j-k} &\equiv &\ (-1)^k {j+k-1 \choose j-1} &\mod p\\
&		&\Leftrightarrow		&\quad \frac{1}{k!} \cdot \prod_{\lambda = 1}^k (p-j-k+\lambda) &\equiv &\ (-1)^k \cdot \frac{1}{k!} \cdot \prod_{\lambda = 1}^k (j+k-\lambda) &\mod p\\
&		&\Leftrightarrow		&\quad \prod_{\lambda = 1}^k (p-j-k+\lambda) &\equiv &\ (-1)^k \prod_{\lambda = 1}^k (j+k-\lambda) &\mod p
\end{align*}
which is clearly true.
\end{proof}

Let $R_n$ be the $(n \times n)$-anti-diagonal matrix over $\F_p$ with ones on the anti-diagonal. The above lemma can be rephrased as
$$\Pasc_p(\F_p)^{-1} = \begin{pmatrix} {p-j \choose p-i} \end{pmatrix}_{i,j} = R_p \cdot \Pasc_p(\F_p)^{\tr} \cdot R_p^{-1}.$$ In order to generalize this result to $(p^r \times p^r)$-matrices, we use the theorem of Lucas.

\begin{thm}[Lucas' Theorem]\label{lucas} Let $p$ be a prime and $a,b \in \N$. Write $a-1 = \sum_{i=0}^r (a_i -1)p^r$ and $b-1 = \sum_{i=0}^r (b_i -1)p^r$ with $1 \leq a_i, b_i \leq p$.
Then
$${a-1 \choose b-1} \equiv {a_r-1 \choose b_r-1} {a_{r-1}-1 \choose b_{r-1}-1} \cdots {a_0-1 \choose b_0-1} \mod p.$$
\end{thm}

Let $p,a,b,a_i,b_i$ be as in Theorem \ref{lucas} and $r \in \N$. By definition, the $(a,b)$ entry of $\Pasc_p(\F_p)^{\otimes r}$ is given by ${a_r-1 \choose b_r-1} {a_{r-1}-1 \choose b_{r-1}-1} \cdots {a_0-1 \choose b_0-1}$ which is, by the theorem, equal to ${a-1 \choose b-1}$ in $\F_p$. This leads to the well known observation 
\begin{equation}\label{powers}
\Pasc_{p^r}(\F_p) = \Pasc_{p}(\F_p)^{\otimes r}
\end{equation}
saying that, modulo a prime, the Pascal triangle has the shape of a Sierpinksi triangle.

\begin{cor}\label{inverse}
Let $p^r$ be a prime power. Then, over $\F_p$, we have
$$\begin{pmatrix} (-1)^{i+j} \cdot {i-1 \choose j-1} \end{pmatrix}_{i,j} = \Pasc_{p^r}(\F_p)^{-1} = \begin{pmatrix} {p^r -i \choose p^r -j} \end{pmatrix}_{i,j}.$$
\end{cor}

\begin{proof}
Using the above results, we find that $\Pasc_{p^r}(\F_p)^{-1}$ is equal to
\begin{align*}
(\Pasc_{p}(\F_p)^{\otimes r})^{-1} &= (\Pasc_{p^r}(\F_p)^{-1})^{\otimes r}\\
&= (R_p \cdot \Pasc_p(\F_p)^{\tr} \cdot R_p^{-1})^{\otimes r} \\
&= R_p^{\otimes r} \cdot (\Pasc_p(\F_p)^{\tr})^{\otimes r} \cdot (R_p^{-1})^{\otimes r}\\
&= R_{p^r} \cdot (\Pasc_p(\F_p)^{\otimes r})^{\tr} \cdot (R_{p^r}^{-1}).
\end{align*}
This implies the desired result.
\end{proof}

\section{Proofs of the Main Results}
\subsection{Proof of Theorem \ref{main1}}
Let $p$ be a prime and $t,n \in \N$. We shall always write $n = mp^r$ with $p\nmid m$. Moreover, let $\F_{p^s} \geq \F_{p^t}$ be a field extension that contains the $m$-th roots of unity. Consider
$$A := \begin{pmatrix} 	0 & 1 & \cdots & 0 \\
			\vdots & \vdots & \ddots & \vdots\\
			0 & 0 & \cdots & 1\\
			1 & 0 & \cdots & 0\end{pmatrix} \in \GL_n(\F_{p^t}).$$
Every circulant matrix over $\F_{p^t}$ can be written as $\left( A^{n-1}v \ A^{n-2}v \ \cdots \ v \right)$ with some (unique) vector $v \in \F_{p^t}^n$. We thus want to determine the number of vectors $v \in \F_{p^t}$ for which $A^{n-1}v, A^{n-2}v, \dots,\ v$ is a basis of $\F_{p^t}^n$.
Let $X \in \GL_n(\F_{p^t})$ and $\widetilde A := XAX^{-1}$. Then we have
\begin{align*}
&A^{n-1}v,\ A^{n-2}, \dots,\ v \mbox{ is a basis of } \F_{p^t}^n\\
\Leftrightarrow \quad &XA^{n-1}v,\ XA^{n-2}v, \dots, \ Xv  \mbox{ is a basis of } \F_{p^t}^n\\
\Leftrightarrow \quad &\widetilde{A}^{n-1}Xv,\ \widetilde{A}^{n-2}Xv, \dots, \ Xv  \mbox{ is a basis of } \F_{p^t}^n.
\end{align*}
Here the last equivalence holds, because $XA^kv = XA^kX^{-1}Xv = \widetilde{A}^kXv$. 
Given an $(n \times n)$-matrix $M$ over $\F_{q}$, with $q$ a power of $p$, let us set
$$ V(M,q) := \{ v \in \F_{q}^n \mid M^{n-1} v,\ M^{n-2}v, \ \dots, \ v  \mbox{ is a basis of } \F_{q}^n\}.$$
We are thus interested in the set $V(A,p^t)$. Observe that we have a bijection 
$$
V(\widetilde{A},p^s) \longrightarrow V(A,p^s), \quad v \longmapsto X^{-1}v.
$$
We may therefore conjugate $A$ into a form $\widetilde{A}$ for which $V(\widetilde{A},p^t)$ is easier to determine. This form $\widetilde{A}$ will actually be the Jordan form of $A$ over the field extension $\F_{p^s} \geq \F_{p^t}$. The above bijection also leads to a bijection
\begin{equation}\label{bij}
V(\widetilde{A},p^s) \cap X \cdot V(A,p^t) \longrightarrow V(A,p^t), \quad v \longmapsto X^{-1}v.
\end{equation}
We are first going to determine the set $V(\widetilde A, p^s)$ and then investigate for which $v \in V(\widetilde A, p^s)$ we have $X^{-1}v \in V(A,p^t)$, that is, we then determine $V(\widetilde{A},p^s) \cap X \cdot V(A,p^t)$.

The characteristic polynomial $\chi_A$ of $A$ is given by $\chi_A = X^n -1$, as one easily finds using Laplace's Theorem. Let $\mu \in \F_{p^s}^*$ be a primitive $m$-th root of unity. Over $\F_{p^s}$, the characteristic polynomial $\chi_A$ of $A$ reads
$$\chi_A \ = \ X^{m p^r} - 1 \ = \  (X^m - 1)^{p^r} \ = \ \prod_{k=1}^{m} (X-\mu^k)^{p^r}.$$
Hence $A$ has the eigenvalues $\mu^k \in \F_{p^s}$, $1 \leq k \leq m$, each with algebraic multiplicity $p^r$. If $v = \begin{pmatrix} v_i \end{pmatrix}_i \in \Eig(A,\mu^k)$, we have
$$
\begin{pmatrix} v_2 \\ \vdots \\ v_n \\ v_1 \end{pmatrix} \ = \ A \cdot v \ = \ \mu^k \cdot v \ = \ \begin{pmatrix} \mu^k \cdot v_1 \\ \vdots \\ \mu^k \cdot v_{n-1} \\ \mu^k \cdot v_n \end{pmatrix} \quad \Rightarrow \quad v = v_1 \cdot \begin{pmatrix} 1 \\ \mu^k \\ \vdots \\ \mu^{(n-1)k}\end{pmatrix}.$$
In particular, $\Eig(A,\mu^k)$ has dimension $1$ and, over $\F_{p^s}$, the Jordan form $\widetilde A$ of $A$ is given by
$$
\widetilde A = \begin{pmatrix} \widetilde A_1 \\ & \ddots \\ & & \widetilde A_{m} \end{pmatrix} \quad \mbox{with} \quad
\widetilde A_b := \begin{pmatrix}
\mu^b & 1\\
      & \mu^b &\ddots\\
      &       & \ddots &1\\
      &       &       & \mu^b
\end{pmatrix} \in \GL_{p^r}(\F_{p^s}).$$

In order to use (\ref{bij}), we have to determine $X\in \GL_n(\F_{p^s})$ such that $XAX^{-1} = \widetilde A$.

For $k \in \N$, let us set 
$$X(k) :=  \begin{pmatrix} (-1)^{i+j} \cdot \mu^ {(i-j)k} \cdot {i-1 \choose j-1} \end{pmatrix} = \begin{pmatrix} \mu^ {(i-j)k} \cdot {p^r - j \choose p^r - i} \end{pmatrix}\in \GL_{p^r}(\F_{p^s}).$$
Here the second equality follows from Corollary \ref{inverse}. We now verify that 
$$X(k)^{-1} = \begin{pmatrix} \mu^ {(i-j)k} \cdot {i-1 \choose j-1} \end{pmatrix}.$$
To this end, note that from (\ref{pasc}) we know that
\begin{equation}\label{sum}\sum_{\lambda=1}^{p^r} (-1)^{\lambda + j} \cdot {i-1 \choose \lambda-1}{\lambda - 1 \choose j - 1} =  \begin{cases}1, \ i = j \\ 0,\ i \not= j  \end{cases}.
\end{equation}
Hence we have
\begin{align*}
& 	& &\begin{pmatrix}\mu^ {(i-j)k} \cdot {i-1 \choose j-1} \end{pmatrix} \cdot \begin{pmatrix} (-1)^{i+j} \cdot \mu^ {(i-j)k} \cdot {i-1 \choose j-1} \end{pmatrix}\\
&= 	& &\begin{pmatrix}\sum_{\lambda=1}^{p^r} \mu^ {(i-\lambda)k} \cdot {i-1 \choose \lambda-1}  (-1)^{\lambda+j} \cdot \mu^ {(\lambda-j)k} \cdot {\lambda-1 \choose j-1} \end{pmatrix}\\
&= 	& &\begin{pmatrix}\mu^ {(i-j)k} \cdot \sum_{\lambda=1}^{p^r} (-1)^{\lambda+j} \cdot {i-1 \choose \lambda-1}{\lambda-1 \choose j-1} \end{pmatrix} \quad \stackrel{(\ref{sum})}{=} \quad I_{p^r}
\end{align*}
where $I_{p^r}$ denotes the $(p^r \times p^r)$-identity matrix.

Let us also set 
$$Y := m^{-1} \cdot \begin{pmatrix} \mu^{-ijp^{r}}\end{pmatrix} \in \GL_m(\F_{p^s}).$$
We verify that the inverse of $Y$ is given by 
$$Y^{-1} = \begin{pmatrix} \mu^{ijp^r}\end{pmatrix}.$$
Indeed, we find
$$
\begin{pmatrix} \mu^{-ijp^r}\end{pmatrix} \cdot \begin{pmatrix} \mu^{ijp^r}\end{pmatrix} \quad 
= \quad \begin{pmatrix} \sum_{\lambda=1}^m \mu^{\lambda(j-i)p^r}\end{pmatrix} \quad
= \quad m\cdot I_m,
$$
since for $i=j$ we have $\mu^{\lambda(j-i)p^r} = 1$ and for $i\not=j$ we observe that $\mu^{(j-i)p^r}$ is a non-trivial $m$-th root of unity and therefore a root of the polynomial $\sum_{\lambda=1}^m X^\lambda$. 

From the matrices $X(k)$ and $Y$ we construct the $(n \times n)$-matrix $X$ by setting
$$X := m^{-1} \begin{pmatrix} \mu^{-abp^r} \cdot X(a)\end{pmatrix}_{1\leq a,b \leq m}.$$
Considering the above computations, one easily verifies that
$$X^{-1} = \begin{pmatrix} \mu^{abp^r} \cdot X(b)^{-1}\end{pmatrix}_{a,b}.$$

We are now going to show that $XAX^{-1} = \widetilde A$ by considering the action of $A$ by multiplication on the columns of $X^{-1}$. Let us write $X^{-1} = \begin{pmatrix} U_1 & \cdots & U_m \end{pmatrix}$ where for $1 \leq b \leq m$ the $(mp^r \times p^r)$-submatrix $U_b$ of $X^{-1}$ is given by 
$$U_b := \begin{pmatrix} \mu^{bp^r} \cdot X(b)^{-1} \\ \vdots \\ \mu^{abp^r} \cdot X(b)^{-1} \\ \vdots \\ X(b)^{-1} \end{pmatrix}.$$
The first column of $\mu^{abp^r} \cdot X(b)^{-1}$ reads
$$\mu^{abp^r} \cdot \begin{pmatrix} 1 \\ \mu^b \\ \mu^{2b} \\ \vdots \\ \mu^{(p^r -1)b} \end{pmatrix}
= \begin{pmatrix} \mu^{ap^r b} \\ \mu^{(ap^r+1)b} \\ \mu^{(ap^r+2)b} \\ \vdots \\ \mu^{(ap^r+(p^r-1))b} \end{pmatrix}
= \begin{pmatrix} \mu^{ap^r b} \\ \mu^{(ap^r+1)b} \\ \mu^{(ap^r+2)b} \\ \vdots \\ \mu^{((a+1)p^r -1))b} \end{pmatrix}.$$
Hence the first column of $U_b$ is given by
$$\begin{pmatrix} \mu^{p^r b} \\ \vdots \\ \mu^{(p^r+i-1)b} \\ \vdots \\ \mu^{(p^r+n-1)b}\end{pmatrix} \in\Eig(A,\mu^b),$$
as desired. In general, the $j$-th column $s_{a,b,j}$ of $\mu^{abp^r} \cdot X(b)^{-1}$ reads
$$s_{a,b,j} = \mu^{abp^r} \cdot \begin{pmatrix} \mu^{(i-j)b}\cdot{i-1 \choose j-1} \end{pmatrix}_{1 \leq i \leq p^r} =  \begin{pmatrix} \mu^{(i-j + ap^r)b}\cdot{i-1 \choose j-1} \end{pmatrix}_{1 \leq i \leq p^r}.$$
For $2 \leq j \leq p^r$ we thus find
\begin{align*}
& & &\ \mu^b \cdot s_{a,b,j} + s_{a,b,j-1}\\
& &=&\ \mu^{abp^r + b} \cdot \begin{pmatrix} \mu^{(i-j)b} \cdot {i-1 \choose j-1}\end{pmatrix}_i + \mu^{abp^r} \cdot \begin{pmatrix} \mu^{(i-(j-1))b} \cdot {i-1 \choose j-2}\end{pmatrix}_i\\
& &=&\ \mu^{(ap^r + 1)b} \cdot \begin{pmatrix} \mu^{(i-j)b} \cdot {i-1 \choose j-1} + \mu^{(i-j)b} \cdot {i-1 \choose j-2}\end{pmatrix}_i\\
& &=&\ \mu^{(ap^r + 1)b} \cdot \begin{pmatrix} \mu^{(i-j)b} \cdot {i \choose j-1}\end{pmatrix}_i.
\end{align*}
Noting that, as a consequence of Lucas' Theorem, ${p^r \choose j-1} = 0 = {0 \choose j-1}$ in $\F_{p^s}$ for $2 \leq j \leq p^r$, we write this as
$$
\mu^b \cdot s_{a,b,j} + s_{a,b,j-1} = 
\begin{pmatrix} \mu^{(2-j+ap^r)b} \cdot {1 \choose j-1} \\ \vdots \\ \mu^{(i+1-j+ap^r)b} \cdot {i \choose j-1} \\ \vdots \\ \mu^{(1-j+(a+1)p^r)b} \cdot {0 \choose j-1}\end{pmatrix}.
$$
Now let $u_{b,j}$ be the $j$-th column of $U_b$, that is $u_{b,j} = \begin{pmatrix} s_{a,b,j} \end{pmatrix}_{1 \leq a \leq m}$.
Then, for $2 \leq j \leq p^r$, we find
$$
\mu^b \cdot u_{b,j} + u_{b,j-1} = 
\begin{pmatrix} \mu^b \cdot s_{1,b,j} + s_{1,b,j-1} \\ \vdots \\ \mu^b \cdot s_{a,b,j} + s_{a,b,j-1}\\ \vdots\\ \mu^b \cdot s_{m,b,j} + s_{m,b,j-1} \end{pmatrix} =
\begin{pmatrix} \mu^{(2-j+p^r)b} \cdot {1 \choose j-1} \\ \vdots \\ \mu^{(1-j+2p^r)b} \cdot {0 \choose j-1}\\ \vdots\\ 
\mu^{(2-j+ap^r)b} \cdot {1 \choose j-1} \\ \vdots \\ \mu^{(1-j+(a+1)p^r)b} \cdot {0 \choose j-1}\\ \vdots \\
\mu^{(2+1-j+mp^r)b} \cdot {1 \choose j-1} \\ \vdots \\ \mu^{(1-j + (m+1)p^r)b} \cdot {0 \choose j-1}\end{pmatrix}
$$
and thus $\mu^b \cdot u_{b,j} + u_{b,j-1} = A \cdot u_{b,j}$.
This proves that $XAX^{-1} = \widetilde A$ indeed has the desired form.

Our next aim is to verify that
\begin{equation}\label{setV}
V(\widetilde A, p^s) = \{ \begin{pmatrix} v_i \end{pmatrix}_i \in \F_{p^s}^n \mid v_{p^r}, v_{2p^r}, \dots, v_{mp^r} \not= 0\}.
\end{equation}
To this end, we need to consider the powers $\widetilde A^k$ for $1 \leq k \leq n-1$.
Obviously, $\widetilde A^k = \diag( \widetilde A_1^k,\dots, \widetilde A_{m}^k)$. We now show that, for $k \in \N$ and $1 \leq b \leq m$, we have
$$\widetilde A_b^k = \begin{pmatrix} \mu^ {(k+i-j)b} \cdot {k \choose j-i}\end{pmatrix}_{i,j}.$$
For $\widetilde A_b^0$ and $\widetilde A_b^1$ this is clear. Suppose that $k \geq 1$ and that the claim is true for $\widetilde A_b^k$. Then we have
\begin{align*}
\widetilde A^{k+1}_b 	
&= \begin{pmatrix} \mu^{(k+i-j)b} \cdot {k \choose j-i} \end{pmatrix}_{i,j} \cdot \begin{pmatrix} \mu^{(1+i-j)b} \cdot {1 \choose j-i} \end{pmatrix}_{i,j}\\
&= \begin{pmatrix} \sum_{\lambda = 1}^{p^r} \mu^{(k+i-\lambda)b} \cdot {k \choose \lambda - i} \cdot  \mu^{(1+ \lambda-j)b} \cdot {1 \choose j-\lambda} \end{pmatrix}_{i,j}\\
&= \begin{pmatrix} \sum_{\lambda = j-1}^{j} \mu^{(k+1+i-j)b} \cdot {k \choose \lambda - i} \cdot {1 \choose j-\lambda} \end{pmatrix}_{i,j}\\
&= \begin{pmatrix} \left( {k \choose j-1-i}+{k \choose j-i} \right) \mu^{(k+1+i-j)b}\end{pmatrix}_{i,j}\\
&= \begin{pmatrix} {k+1 \choose j-i}  \mu^{(k+1+i-j)b}\end{pmatrix}_{i,j}
\end{align*}
as claimed.

Recall that for $w = \begin{pmatrix} w_i \end{pmatrix} \in \F_{p^s}^{p^r}$, we write $T(w) =  \begin{pmatrix} w_{p^r-j+i}\end{pmatrix}_{1\leq i,j \leq p^r}$ where $w_\lambda = 0$ for $\lambda \leq 0$. Observe that $T(w)$ is regular, if and only if $w_{p^r} \not= 0$. Using the matrix $T(w)$, we can write
\begin{align*}
\widetilde A^k_b \cdot w &= \begin{pmatrix} \sum_{\lambda=1}^{p^r} \mu^{(k+i-\lambda)b}\cdot{k \choose \lambda - i}\cdot w_\lambda \end{pmatrix}\\
&= \begin{pmatrix} \sum_{\lambda=i}^{k+i} \mu^{(k+i-\lambda)b}\cdot{k \choose \lambda - i}\cdot w_\lambda \end{pmatrix}\\
&= \begin{pmatrix} \sum_{\lambda=p^r-k}^{p^r} \mu^{(k-p^r+\lambda)b}\cdot{k \choose p^r - \lambda}\cdot w_{p^r - \lambda +i} \end{pmatrix}\\
&= \begin{pmatrix} \sum_{\lambda=1}^{p^r} \mu^{(k-p^r+\lambda)b}\cdot{k \choose p^r - \lambda}\cdot w_{p^r - \lambda + i} \end{pmatrix}\\
&= T(w) \cdot \begin{pmatrix} \mu^{(k- p^r +i)b}{k \choose p^r-i}\end{pmatrix}.
\end{align*}
We thus obtain
\begin{align*}
\begin{pmatrix} \widetilde A_b^{p^r-1}w & \widetilde A_b^{p^r-2}w & \cdots & w \end{pmatrix} &= T(w) \cdot \begin{pmatrix} \mu^{(p^r-j-p^r+i)b}{p^r-j \choose p^r -i}  \end{pmatrix}\\
&= T(w) \cdot X(b)^{-1}.
\end{align*}
We note that $\widetilde A_b^{kp^r} = \mu^{kbp^r} \cdot I_{p^r}$ so that $\widetilde A_b^{kp^r+l} = \mu^{kp^rb} \cdot \widetilde A_b^l$. Hence
\begin{align*}
 & \begin{pmatrix} \widetilde A_b^{n-1} w & \widetilde A_b^{n-2} w & \cdots & w \end{pmatrix} = \begin{pmatrix} \widetilde A_b^{mp^r-1} w & \widetilde A_b^{mp^r-2} w & \cdots & w \end{pmatrix}\\
=&\ T(w) \cdot \begin{pmatrix} \mu^{(m-1)p^r b} X(b) & \mu^{(m-2)p^r b} X(b) & \cdots & X(b) \end{pmatrix}.
\end{align*}
Now let $v = \begin{pmatrix} v_{(i)}\end{pmatrix}_{1\leq i \leq m}  \in \F_{p^s}^{n}$ where $v_{(i)} = \begin{pmatrix} v_{i,j}\end{pmatrix}_{1 \leq j \leq p^r} \in \F_{p^s}^{p^r}$. Then $\widetilde A^k v = \begin{pmatrix} \widetilde A_i^k v_{(i)} \end{pmatrix}_{1 \leq i \leq m}$ so that $\begin{pmatrix} \widetilde A^{n-1} v & \widetilde A^{n-2} v & \cdots & v\end{pmatrix}$ reads
$$\begin{pmatrix}
 	\widetilde A_1^{n-1} v_{(1)} & \widetilde A_1^{n-2} v_{(1)} & \cdots & v_{(1)}\\
	\widetilde A_2^{n-1} v_{(2)} & \widetilde A_2^{n-2} v_{(2)} & \cdots & v_{(2)}\\
	\vdots& \vdots & \ddots & \vdots \\
	\widetilde A_m^{n-1} v_{(m)} & \widetilde A_m^{n-2} v_{(m)} & \cdots & v_{(m)}\end{pmatrix}.$$
This can be written as the product
\begin{align*}
&\begin{pmatrix} T(v_{(1)})\\ & T(v_{(2)})\\ & &\ddots\\ & & & T(v_{(m)}) \end{pmatrix}\\
&\ \cdot
	\begin{pmatrix}
 	\mu^{(m-1)p^r} X(1) & \mu^{(m-2)p^r } X(1) & \cdots & X(1) \\
	\mu^{2(m-1)p^r} X(2) & \mu^{2(m-2)p^r } X(2) & \cdots & X(2)\\
	\vdots & \vdots & \ddots & \vdots\\
	X(m) & X(m) & \cdots & X(m)
	\end{pmatrix}
\end{align*}
where the right matrix is exactly the matrix $m \cdot X$. Hence
\begin{equation}\label{product}
\begin{pmatrix} \widetilde A^{n-1} v & \widetilde A^{n-2} v & \cdots & v\end{pmatrix} = \diag(T(v_{(1)}), \dots, T(v_{(m)})) \cdot m \cdot X.
\end{equation}
Clearly, this product is regular, if and only if so  are $T(v_{(1)}),\dots,T(v_{(m)})$, that is, if and only if $v_{1,p^r},\dots,v_{m,p^r} \not= 0$. This proves (\ref{setV}).

Now we have to investigate for which $v \in \F_{p^s}^n$ we have $X^{-1} v \in \F_{p^t}^n$. Recall that we write $F^t: \F_{p^s} \rightarrow \F_{p^s}$, $x \mapsto x^{p^t}$ for the relative Frobenius Homomorphism. We use the standard result that $\F_{p^t} = \{ x \in \F_{p^s} \mid F^t(x) = x\}$ to find
\begin{align*}
X^{-1} v \in \F_{p^t}^n &\Leftrightarrow F^t(X^{-1} v) = X^{-1} v \\
&\Leftrightarrow F^t(X)^{-1} F^t(v) = X^{-1} v\\
&\Leftrightarrow F^t(X) X^{-1} v = F^t(v).
\end{align*}
Consider the matrix $F^t(X)X^{-1}$. One easily sees that 
$$F^t(X) = m^{-1}\begin{pmatrix} \mu^{-(a p^t) bp^r} \cdot X(ap^t)\end{pmatrix}$$
so that $F^t(X) = PX$ where $P \in \GL_n(\F_p)$ is a permutation matrix consisting of $(p^r \times p^r)$-blocks which describes the permutation $\sigma: \Z/m\Z \rightarrow \Z/m\Z$, $a \mapsto p^{-t} a$. Observe that the cycle structure of $\sigma^{-1}$ and hence also of $\sigma$ is described in Proposition \ref{cycle}. We have
\begin{align*}
F^t(X) X^{-1} v = F^t(v) &\Leftrightarrow Pv = F^t(v) \\
&\Leftrightarrow v_{(\sigma(i))} = F^t(v_{(i)}), \ 1 \leq i \leq m.
\end{align*}
Let us write $\sigma = \sigma_1 \cdots \sigma_l$ with disjoint cycles $\sigma_k$ of length $m_k$. Clearly
$$v_{(\sigma(i))} = F^t(v_{(i)}) \quad \Leftrightarrow \quad v_{(\sigma_k(i))} = F^t(v_{(i)}),\ 1 \leq k \leq l.$$
Writing $\sigma_k=(s_{k,1}, \dots, s_{k,m_k})$, the condition $v_{(\sigma_k(i))} = F^t(v_{(i)})$, $1 \leq i \leq m$ means 
$$v_{(s_{k,i})} = F^t(v_{(s_{k,{i-1}})}),\ 2 \leq i \leq m_k \ \mbox{and} \ v_{(s_{k,1})} = F^t(v_{(s_k,{m_k})}).$$
By a simple substitution, we can rewrite this as
$$v_{(s_{k,i})} = F^{(i-1)t}(v_{(s_{k,1})}),\ 2 \leq i \leq m_k \ \mbox{and} \ v_{(s_{k,1})} = F^{t m_k}(v_{(s_k,1)})$$
which in turn can be written as
$$v_{(s_{k,i})} = v_{(s_{k,1})}^{p^{(i-1)t}},\ 2 \leq i \leq m_k \ \mbox{and} \ v_{(s_{k,1})} \in \F_{p^{t m_k}}^{p^r}.$$
Hence the vectors $v = \begin{pmatrix} v_{(i)}\end{pmatrix} \in V(\widetilde A, p^s)$ for which $X^{-1}v \in \F_{p^t}^n$ are precisely the ones that satisfy
\begin{align}
\label{cond1}&v_{1,p^r}, v_{2,p^r}, \dots v_{m,p^r} \not= 0,\\
\label{cond2}& v_{(s_{k,1})} \in \F_{p^{t m_k}}^{p^r},\ 1 \leq k \leq l,\\
\label{cond3}& v_{(s_{k,i})} = F^{(i-1)t}(v_{(s_{k,1})}),\ 2 \leq i \leq m_k,\ 1 \leq k \leq l.
\end{align} 
One easily verifies that there are exactly 
$$\prod_{k=1}^l (p^{t m_k} - 1)(p^{t m_k})^{p^r-1} = \prod_{k=1}^l(p^{t m_k}-1)(p^{t m_k (p^r-1)})$$
such vectors. Using Proposition \ref{cycle}, which gives the cycle structure of $\sigma$, i.e. the numbers $m_k$,  we thus obtain part(i) of Theorem \ref{main1}.

Using our results so far, it is not difficult to proof part (ii). The group $\RC_n(\F_{p^t})$ consists exactly of the matrices
$$\begin{pmatrix} A^{n-1}X^{-1}v & A^{n-2}X^{-1}v & \cdots & X^{-1}v \end{pmatrix}$$
with $v \in \F_{p^s}^{n}$ satisfying (\ref{cond1}), (\ref{cond2}), (\ref{cond3}). Observe that
\begin{align*}
&\ X \begin{pmatrix} A^{n-1}X^{-1}v & A^{n-2}X^{-1}v & \cdots & X^{-1}v \end{pmatrix} X^{-1}\\
= &\ \begin{pmatrix} \widetilde A^{n-1} v & \widetilde A^{n-2} v & \cdots & v \end{pmatrix} X^{-1}\\
= &\ m \cdot \diag(T(v_{(1)}), \dots, T(v_{(m)})), \mbox{ by (\ref{product})}.
\end{align*}
Since $m \in \F_{p^t}^*$, this proves part (ii) of Theorem \ref{main1}.

\subsection{Proof of Theorem \ref{main2}}
Let $a,n \in \N$ and consider the group $\RC_n(\Z/a\Z)$ of regular circulant $(n \times n)$-matrices over $\Z/a\Z$. Let $a = p_1^{t_1} \cdots p_k^{t_k}$ be the prime factorization of $a$. Then the isomorphism
$$\GL_n(\Z /a \Z) \stackrel{\cong}{\longrightarrow} \GL_n(\Z / p_1^{t_1} \Z) \times \cdots \times \GL_n(\Z / p_k^{t_k} \Z)$$
leads to an isomorphism
\begin{equation}\label{iso}
\RC_n(\Z /a \Z) \stackrel{\cong}{\longrightarrow} \RC_n(\Z / p_1^{t_1} \Z) \times \cdots \times \RC_n(\Z / p_k^{t_k} \Z).
\end{equation}
Hence we only have to consider $\RC_n(\Z/p^t \Z)$ with $t \in \N$ and $p$ prime.

We have an exact sequence
\begin{equation}\label{seq}
1 \longrightarrow K_t \longrightarrow \GL_n(\Z / p^{t+1} \Z) \longrightarrow \GL_n(\Z / p^t \Z) \longrightarrow 1
\end{equation}
where $K_t = \{ I_n + p^t \begin{pmatrix}a_{ij}\end{pmatrix} \mid a_{ij} \in \Z / p^{t+1} \Z\}$. Every element of $K_t$ can be written in a unique way as $I_n + p^t \begin{pmatrix} a_{ij} \end{pmatrix}$ with $a_{ij} \in \{ 0,1,\dots,p-1 \}$. Morever, $I_n + p^t \begin{pmatrix} a_{ij} \end{pmatrix}$ is circulant, if and only if $\begin{pmatrix} a_{ij}\end{pmatrix}$ is. Hence we find
$$|K_t \cap \RC_n(\Z / p^{t+1} \Z)| = p^n.$$

Now let $C \in \RC_n(\Z / a \Z)$. We can easily lift $C$ to $\RC_n(\Z/p^{t+1} \Z)$ via the natural projection as follows. Let $\begin{pmatrix} c_i \end{pmatrix}_{1\leq i \leq n}$ be the first column of $C$. For every $c_i \in \Z / p^t \Z$, let $\tilde c_i \in \Z/ p^{t+1}\Z$ be a lift and let $\widetilde C$ be the circulant matrix over $\Z / p^{t+1}\Z$ with first column $\begin{pmatrix} \tilde c_i \end{pmatrix}$. Obviously, $\widetilde C$ is a lift of $C$. Moreover, $\det(\widetilde C) \equiv \det(C) \not\equiv 0 \mod p$ so that $\det(\widetilde C)$ is a unit and $\widetilde C$ is regular. Hence $\widetilde C \in \RC_n(\Z / p^{t+1}\Z)$. It follows that the sequence (\ref{seq}) leads to an exact sequence
$$1 \longrightarrow K_t \cap \RC_n(\Z/ p^{t+1}\Z) \longrightarrow \RC_n(\Z / p^{t+1} \Z) \longrightarrow \RC_n(\Z / p^t \Z) \longrightarrow 1$$
so that
\begin{align*}
|\RC_n(\Z / p^{t+1}\Z)| &= |K_t \cap \RC_n(\Z/ p^{t+1}\Z)| \cdot |\RC_n(\Z / p^t \Z)|\\
&= p^n \cdot |\RC_n(\Z / p^t \Z)|.
\end{align*}
Writing $n = p^r m$ such that $p \nmid m$, we know by Theorem \ref{main1} that
$$|\RC_n(\Z / p \Z)| = \prod_{d \mid m} (p^{\ord_d(p)}-1)^{\phi(d) / \ord_d(p)}  \cdot p^{(p^r-1) \phi(d)}.$$
By induction we obtain
$$|\RC_n(\Z / p^t \Z)| = p^{(t-1)n} \cdot \prod_{d \mid m} (p^{\ord_d(p)}-1)^{\phi(d) / \ord_d(p)}  \cdot p^{(p^r-1) \phi(d)}.$$
Finally, by isomorphism (\ref{iso}),  this leads to Theorem \ref{main2}.

\end{document}